\documentclass{elsart}
\usepackage{amsmath}
\usepackage{amssymb}
\usepackage{graphicx,subfigure}
\usepackage{bm}

\newcommand{\Hil}[0]{
\mathcal{H} 
}
\newcommand{\norm}[2]{
\left\| #2 \right\|_{#1}
}
\newcommand{\LtRd}[0]{
{L^2 \left( {\mathbb R}^d \right)}
}
\newcommand{\HS}[0]{
{\mathcal HS}
}
\newcommand{\MM}[0]{
{\bf M} 
}

\newcommand{\RR}[0]{
\mathbb{R} 
}
\newcommand{\Rd}[0]{
{\mathbb{R}^d}
}
\newcommand{\BL}[0]{
{\mathcal B}
}
\newcommand{\LtRtd}[0]{
{L^2 \left( {\mathbb R}^{2d} \right)}
}
\newcommand{\CC}[0]{
{\mathbb{C}}
}
\newcommand{\BB}[0]{
{\mathfrak B}
}

\newtheorem{theorem}{Theorem}[section]
\newtheorem{definition}{Definition}[section]
\newtheorem{proposition}[theorem]{Proposition}
\newtheorem{lemma}[theorem]{Lemma}
\newtheorem{corollary}[theorem]{Corollary}
\newtheorem{conj.}[theorem]{Conjecture}
\newtheorem{Bsp.}{Example}[section]
\newcommand{\Bsp}[1]{
\begin{Bsp.} \label{#1} \bf : \end{Bsp.}
}
\newcommand{\Rem}[0]{

\hspace{-\parindent}%
{\rm \bf%
Remark:}
}
\newenvironment{proof}{\noindent \bf Proof: \rm}{$ \hspace{\stretch{1}} \Box $

\vspace{5mm}}
\begin{document}
\begin{frontmatter}
\title{Matrix Representation of Operators Using Frames
}
\author{Peter Balazs}
\address{Austrian Academy of Sciences, Acoustic Research Institute,
         Reichsratsstrasse 17,A-1010 Vienna, Austria}
\date{\today}

\begin{abstract}
In this paper it is investigated how to find a matrix representation of operators on a Hilbert space $\Hil$ with Bessel sequences, frames and Riesz bases. In many applications these sequences are often preferable to orthonormal bases (ONBs). Therefore it is useful to extend the known method of matrix representation by using these sequences instead of ONBs for these application areas. 
We will give basic definitions of the functions connecting infinite matrices defining bounded operators on $l^2$ and operators on $\Hil$. We will show some structural results and give some examples. 
Furthermore in the case of Riesz bases we prove that those functions are 
isomorphisms. Finally we are going to apply this idea to the connection of Hilbert-Schmidt operators and Frobenius matrices.
\end{abstract}
\vspace{3mm}
\end{frontmatter}

\section{Introduction}

The relevance of signal processing in today's life is clearly evident.
Without exaggeration it can be said, that any advance in signal processing
sciences directly leads to an application in technology and information
processing. Without signal processing methods several modern
technologies would not be possible, like mobile phone, UMTS, xDSL
or digital television. 

The mathematical background for today's signal processing applications are \em Gabor \em \cite{feistro1} , \em wavelet \em \cite{daubech1} and \em sampling theory \em \cite{befe01}. 
A signal is sampled and then analyzed using a Gabor respectively wavelet system. 
Many applications use a modification on the coefficients obtained from the analysis operation \cite{Kronland91application,hlawatgabfilt1}. 
For them not only an analysis but also a synthesis operation is needed. If the coefficients are not changed, the result of the synthesis should be the original signal, therefore so-called \em perfect reconstruction \em is needed. One way to achieve that is to analyze the signal using \em orthonormal bases \em (ONBs). In this case the analysis of a function is simply the correlation of the signal $f$ with each basis element $e_k$,  $f \mapsto (c_k) := \left( \left< f , e_k \right> \right)$ . The synthesis that gives perfect reconstruction
is simply the (possibly infinite) linear combination of the basis elements using the coefficients $c = (c_k)$, $c \mapsto \sum \limits_k c_k e_k$.

From practical experience it soon has become 
apparent that the concept of an orthonormal basis is not always useful. Sometimes it is more important for a decomposing set to have other special properties rather than guaranteeing unique coefficients. For example it is impossible to have good time-frequency localization for Gabor ONBs or a wavelet ONB with a mother wavelet which has exponentially decay and is infinitely often differentiable with bounded derivatives \cite{ole1}.
Furthermore suitable ONBs are often difficult to construct in a numerical efficient way. 
This led to the concept of frames, which was introduced by Duffin and Schaefer \cite{duffschaef1}. It was made popular by Daubechies \cite{daubech1}, and today it is one of the most important foundations of Gabor \cite{feistro1}, wavelet \cite{aliant1} and sampling theory \cite{aldrgroech1}. In signal processing applications frames have received more and more attention \cite{boelc1,vettkov1}. 

Models in physics \cite{aliant1} and other application areas, for example in sound vibration analysis \cite{kreizxxl1}, are mostly continuous models. A lot of problems there can be formulated as operator theory problems, for example in differential or integral equations. To be able to work numerically the operators have to be discretized. One way to do this is to find (possibly infinite) matrices describing these operators using ONBs. In this paper we will investigate a way to describe an operator as a 
 matrix using frames. This kind of 'sampling of operators' (compare to \cite{goetzwal1}) is especially important for application areas, where frames are heavily used, so that the link between model and discretization is kept.
For implementations operator equations can be transformed in a finite, discrete problem with the finite section method \cite{gohberg1} in the same way as in the ONB case.

The standard matrix description \cite{conw1} of operators $O$ using an ONB $(e_k)$ is by constructing an matrix $M$ with the entries $M_{j,k} = \left<O e_k, e_j\right>$. In \cite{olepinv} a concept was presented, where an operator $R$ is described by the matrix $\left( \left< R \phi_j , \tilde{\phi}_i \right> \right)_{i,j}$ with $(\phi_i)$ being a frame and $(\tilde{\phi}_i)$ its canonical dual. Recently such a kind of representation is used for the description of operators in \cite{GroechSjoe1} using Gabor frames and \cite{strohmpdo1} using linear independent Gabor systems. In this paper we are going to develop this idea in full generality for Bessel sequences, frames and Riesz sequences and also look at the dual function which assigns an operator to a matrix. 

This paper is organized as follows: In Section \ref{sec:notprelim0} we collect results and notation we need. 
Section \ref{sec:descropfram0} gives the basic definitions and properties for Bessel sequences and frames. Matrix representation with Riesz bases is covered in Section \ref{sec:matrepriesz0}. In Section \ref{sec:descrhsopfram0} the connection of Frobenius matrices and Hilbert-Schmidt is investigated. Section \ref{sec:concpersp0} finishes the paper with perspectives.

\section{Notation and Preliminaries} \label{sec:notprelim0}

\subsection{Hilbert spaces and Operators }

We will give only a short review, for details refer to \cite{conw1}. We will denote infinite dimensional Hilbert spaces by $\Hil$ 
and their inner product with $<.,.>$, which is linear in the first coordinate. 
Let 
$\BL(\Hil_1,\Hil_2)$ denote the set of all linear and bounded operators from $\Hil_1$ to $\Hil_2$. 
With the \em operator norm\em, $ \norm{Op}{A} = \sup \limits_{\norm{\Hil_1}{x} \le 1} \left\{ \norm{\Hil_2}{A (x)} \right\} $, this set is a Banach space. We will denote the composition of two operators $A: \Hil_1 \rightarrow \Hil_2$ and $B: \Hil_2 \rightarrow \Hil_3$ by $B \circ A : \Hil_1 \rightarrow \Hil_3$ and the adjoint of an operator $A$ by $A^*$, so that $\left< A x, y\right> = \left< x, A^* y \right>$ for all $x,y \in \Hil$. 

Furthermore we will denote the range of an operator $A$ by $ran(O)$
 and its kernel by $ker(A)$. 
An example for a Hilbert space is the sequence space $l^2$ consisting of all square-summable sequences in $\CC$
with the inner product $\left< c , d \right> = \sum \limits_{k} c_k \cdot \overline{d_k}$. 
We will use the canonical basis $\Delta = (\delta_k)$ for sequence spaces
, where $\left( \delta_k \right)_n = \delta_{k,n}$, using the \em Kronecker symbol\em :
$\delta_{k,n} = 
\left\{ \begin{array}{c c} 1 & k = n \\ 0 & otherwise \end{array} \right.$.

Remember that a linear function between Banach algebras $\varphi : \BB_1 \rightarrow \BB_2$ is called a \em Banach algebra homomorphism\em , if it is also \em multiplicative\em , i.e. for all $x,y \in \BB_1$ we have $\varphi(x \cdot y ) = \varphi(x) \cdot \varphi(y)$.
It is called a \em monomorphism\em , if it is also injective.

\begin{definition} \label{sec:kronprod1} Let $X,Y,Z$ be sets, $f : X \rightarrow Z$, $g : Y \rightarrow Z $ be arbitrary functions. The \em Kronecker product \em $ \otimes_{o} : X \times Y \rightarrow Z$ is defined by 
$$ \left( f \otimes_{o} g \right) (x,y) = f(x) \cdot g(y) .$$

Let $f \in \Hil_1$, $g \in \Hil_2$ then define the \em inner tensor product \em as an operator from $\Hil_2$ to $\Hil_1$ by
$$
\left( f \otimes_{i} \overline{g} \right) (h) = \left< h, g \right> f  \mbox{ for } h \in \Hil_2.$$
\end{definition}
We will often write $f \otimes  g$ instead of $ f \otimes_{o} g$ or $ f \otimes_{i} g$ , if there is no chance of misinterpretation. 

\subsubsection{Hilbert Schmidt Operators} \label{sec:hilbertsch}

A bounded operator $T \in \BL (\Hil_1,\Hil_2)$ is called a \em Hilbert-Schmidt \em ($\HS$) operator  if there exists an ONB $( e_n ) \subseteq \Hil_1$ such that \index{operator!Hilbert Schmidt} 
$$ \norm{\HS}{T} := \sqrt{ \sum \limits_{n=1}^{\infty} \norm{\Hil_2}{T e_n}^2} < \infty $$
Let ${\mathcal HS}(\Hil_1, \Hil_2)$ denote the space of Hilbert Schmidt operators from $\Hil_1$ to $\Hil_2$.

This definition is independent of the choice of the ONB.  The class of Hilbert-Schmidt operators is a Hilbert space of the compact operators with the following properties: 
\begin{itemize}
\item $ \norm{Op}{T} \le \norm{\HS}{T} $
\item $\norm{\HS}{T} = \norm{\HS}{T^*}$, and $T \in \HS$ $\Longleftrightarrow$ $T^* \in \HS$.
\item 
If $T \in \HS$ and $A \in \BL$, then $T A$ and $A T \in \HS$. $\norm{\HS}{A T} \le \norm{Op}{A} \norm{\HS}{T}$ and $\norm{\HS}{T A} \le \norm{Op}{A} \norm{\HS}{T}$. 
\end{itemize}

For more details on this class of compact operators refer to \cite{schatt1} or \cite{wern1}.

\subsection{Frames}

For more details and proofs for this section refer e.g. to \cite{Casaz1,ole1,daubech1,Groech1}.

A sequence ${\Psi} = \left( \psi_k | k \in K \right)$ 
 is called a \em frame \em for the Hilbert space $\Hil$, if constants $A,B > 0$ exist, such that 
\begin{equation} \label{sec:framprop1} A \cdot \norm{\Hil}{f}^2 \le \sum \limits_k \left| \left< f, \psi_k \right> \right|^2 \le B \cdot  \norm{\Hil}{f}^2  \ \forall \ f \in \Hil
\end{equation} 
Here $A$ is called a \em lower\em , $B$ an \em upper frame bound\em . 
If the bounds can be chosen such that $A=B$ the frame is called \em tight\em .

A sequence $\Psi = (\psi_k)$ is called a \em Bessel sequence \em with Bessel bound $B$ if it fulfills the right inequality above:
\begin{equation} \label{sec:bessel1} \sum \limits_k \left| \left< f, \psi_k \right> \right|^2 \le B \cdot  \norm{\Hil}{f}^2  \ \forall \ f \in \Hil
\end{equation}

The index set will be omitted in the following, if no distinction is necessary. 

For a Bessel sequence, $\Psi = ( \psi_k )$, let $C_{\Psi} : \Hil \rightarrow l^2 ( K )$ be the \em analysis  operator \em
$ C_{\Psi} ( f ) = \left( \left< f , \psi_k \right> \right)_k$. 
Let $D_{\Psi} : l^2( K ) \rightarrow \Hil $ be the  \em synthesis operator \em
$ D_{\Psi} \left( \left( c_k \right) \right) = \sum \limits_k c_k \cdot \psi_k $. 
Let $S_{\Psi} : \Hil  \rightarrow \Hil $ be the \em (associated) frame  operator \em
$ S_{\Psi} ( f  ) = \sum \limits_k  \left< f , \psi_k \right> \cdot \psi_k $. 
To simplify notation we will just write $S$ for $S_{\Psi}$, $C$ for $C_{\Psi}$ and $D$ for $D_{\Psi}$, if it is not necessary to distinguish different frames. We will use the notation $S_{\Psi,\Phi} = D_\Psi \circ C_\Phi$.
$C$ and $D$ are adjoint to each other, $D = C^*$ with $\norm{Op}{D} = \norm{Op}{C} \le \sqrt{B}$.  The series $\sum \limits_k c_k \cdot \psi_k$ converges unconditionally for all $(c_k) \in l^2$. 

For a frame ${\Psi} = ( \psi_k )$ with bounds $A,B$, $C$ is a bounded, injective operator with closed range and $S = C^*C = DD^*$ is a positive invertible operator satisfying $A I_\Hil \le S \le B I_\Hil$ and $B^{-1} I_\Hil \le S^{-1} \le A^{-1} I_\Hil$. Even more we can find an expansion for every member of $\Hil$: 
The sequence $\tilde{\Psi} = \left( \tilde{\psi}_k \right) = \left( S^{-1} \psi_k \right)$ 
is a frame with frame bounds $B^{-1}$, $A^{-1} > 0$, the so called \em canonical dual frame\em .
 Every $f \in \Hil$ has the 
 expansions
$ f = \sum \limits_{k \in K} \left< f, \tilde{\psi}_k \right> \psi_k $
and 
$ f = \sum \limits_{k \in K} \left< f, \psi_k \right> \tilde{\psi}_k $
where both sums converge unconditionally in $\Hil$.

Remember that a
sequence $( e_k )$ is called a \em (Schauder) basis \em 
for $\Hil$, if for all $f \in \Hil$ there are unique coefficients $(c_k)$ such that
$ f = \sum \limits_k c_k \phi_k $. Also
two sequences $(\psi_k)$, $(\phi_k)$ are called \em biorthogonal \em 
 if 
$\left< \psi_k, \phi_j\right> = \delta_{kj}$ for all $h,j$.

A complete sequence $( \psi_k)$ in $\Hil$  is called a \em Riesz basis \em if 
there exist constants $A$, $B >0$ such that the inequalities
$$ A \norm{2}{c}^2 \le \norm{\Hil}{\sum \limits_{k \in K} c_k \psi_k}^2 \le B \norm{2}{c}^2 $$
hold for all finite sequences $(c_k)$.

For a frame $( \psi_k )$ the following conditions are equivalent:
$(i)$
$( \psi_k )$ is a Riesz basis for $\Hil$.
$(ii)$ The coefficients $( c_k ) \in l^2$ for the series expansion with $( \psi_k )$ are unique. So the synthesis operator $D$ is injective.
$(iii)$ The analysis operator $C$ is surjective.
$(iv)$ $( \psi_k )$ and $( \tilde{\psi}_k )$ are biorthogonal.

Let $\Psi = ( \psi_k )$ and $\Phi = ( \phi_k)$ be two sequences in $\Hil$. 
The \em Gram matrix \em$G_{\Psi, \Phi}$ for these sequences is given by $\left( G_{\Psi, \Phi} \right)_{j,m} = \left< \phi_m , \psi_j \right>$, $j,m \in K$. 
We denote $G_{\Psi,\Psi}$ by $G_\Psi$.
We can look at the operator induced by the Gram matrix, defined for $c \in l^2$ formally as
$ ( G_{\Psi , \Phi} c )_j = \sum \limits_k c_k \left< \phi_k , \psi_j \right> $. 
Clearly for two Bessel sequences it is well defined as linear bounded operator, because
$$ ( G_{\Psi , \Phi} c )_j = \sum \limits_k c_k \left< \phi_k , \psi_j \right> =´ \left< \sum \limits_k c_k  \phi_k , \psi_j \right> =   \left( \left(  C_{\Psi} \circ D_{\Phi} \right) c \right)_j $$
and therefore  
$\norm{Op}{G_{\Psi , \Phi}} \le \norm{Op}{C_{\Psi}} \norm{Op}{D_{\Phi}} \le  B$. A frame is a Riesz sequence if and only if the Gram matrix defines a bounded and invertible operator on $l^2$.

\section{Representing Operators with Frames} \label{sec:rowf0}

Let $(\psi_k)$ be a frame in $\Hil_1$. An existing operator $U \in \BL(\Hil_1,\Hil_2)$ is uniquely determined by its images of the frame elements. For $f = \sum \limits_k c_k \psi_k$ 
$$U (f) = U ( \sum_k c_k \psi_k ) =  \sum_k c_k U (\psi_k).$$ 

On the other hand, contrary to the case for ONBs, we cannot just choose a Bessel sequence $(\eta_k)$ and define an operator just by choosing $V(\psi_k) := \eta_k$ and setting $V( \sum \limits_k c_k \psi_k ) = \sum \limits_k c_k \eta_k $. This is in general not well-defined. Only if 
$$ \sum \limits_k c_k \psi_k =  \sum \limits_k d_k \psi_k \Longrightarrow \sum \limits_k c_k \eta_k = \sum \limits_k d_k \eta_k$$
this definition is non-ambiguous, i.e. if $ker \left( D_{\psi_k} \right) \subseteq ker \left( D_{\eta_k} \right)$. This condition is certainly fulfilled, if $D_{\psi_k}$ is injective, i.e. for Riesz bases. 

This problem can be avoided by using the following definition  
\begin{equation} \label{sec:defopbyfram1} V (f) :=  \sum_k \left< f , \tilde{\psi}_k \right> \eta_k.
\end{equation}
As $(\eta_k)$ forms a Bessel sequence, the right hand side of Eq. (\ref{sec:defopbyfram1}) is well-defined. It is clearly linear, and it is bounded.  
The Bessel condition, Eq. \ref{sec:bessel1}, is necessary in the case of ONBs to get a bounded operator, too \cite{conw1}. 
But contrary to the ONB case, here, in general, $V(\psi_k) \not= \eta_k$.

Instead of changing the sequence with which the coefficients are resynthezised, an operator can also be described by changing the coefficients, as presented in the following sections.

\subsection{Matrix Representation} \label{sec:descropfram0}

For orthonormal sequence it is well known, that operators can be uniquely described by a matrix representation \cite{gohberg1}. The same can be constructed with frames and their duals. 
Recall the definition of the operator defined by a (possibly infinite) matrix: $\left( M c\right)_j = \sum \limits_k M_{j,k} c_k$.
We will start with the more general case of Bessel sequences. Note that we will use the notation $\norm{\Hil_1 \rightarrow \Hil_2}{.}$ for the operator norm in $\BL(\Hil_1, \Hil_2)$ to be able to distinguish between different operator norms.
\begin{theorem} \label{sec:matbyfram1} Let $\Psi = (\psi_k)$ be a Bessel sequence in $\Hil_1$ with bound $B$, $\Phi = (\phi_k)$ in $\Hil_2$ with $B'$.
\begin{enumerate} \item Let $O : \Hil_1 \rightarrow \Hil_2$ be a bounded, linear operator. Then the infinite matrix 
$$ 
{\left( {\mathcal M}^{(\Phi , \Psi)} \left( O \right) \right)}_{m,n} = 
\left<O \psi_n, \phi_m \right>$$%
defines a bounded operator from $l^2$ to $l^2$ with $\norm{l^2 \rightarrow l^2}{\mathcal M} \le \sqrt{B \cdot B'} \cdot \norm{\Hil_1 \rightarrow \Hil_2}{O}$.
As an operator $l^2 \rightarrow l^2$ 
$$ {\mathcal M}^{(\Phi , \Psi)} \left( O \right) = C_{\Phi} \circ O \circ D_{\Psi} $$
This means the function ${\mathcal M}^{(\Phi , \Psi)} : \BL(\Hil_1,\Hil_2) \rightarrow \BL(l^2,l^2)$ is a well-defined bounded operator.
\item On the other hand let $M$ be an infinite matrix defining a bounded operator from $l^2$ to $l^2$, $\left(M c\right)_i = \sum \limits_k M_{i,k} c_k$. Then the operator $\mathcal{O}^{(\Phi , \Psi)}$ defined by 
$$ \left( \mathcal{O}^{(\Phi , \Psi)} \left( M \right)\right) h = \sum \limits_k  \left( \sum \limits_j M_{k,j} \left<h, \psi_j\right> \right) \phi_k \mbox{, for } h \in \Hil_1$$  
is a bounded operator from $\Hil_1$ to $\Hil_2$ with 
$$\norm{\Hil_1 \rightarrow \Hil_2}{\mathcal{O}^{(\Phi , \Psi)} \left( M \right)} \le \sqrt{B \cdot B'} \norm{l^2 \rightarrow l^2}{M}.$$
$$ \mathcal{O}^{(\Phi , \Psi)} (M) = D_{\Phi} \circ M \circ C_{\Psi} = \sum \limits_k  \sum \limits_j M_{k,j} \cdot \phi_k \otimes_i \overline{\psi}_j $$
This means the function ${\mathcal O}^{(\Phi , \Psi)} : \BL(l^2,l^2) \rightarrow \BL(\Hil_1,\Hil_2)$ is a well-defined bounded operator.
\end{enumerate} 
\end{theorem}
\begin{proof} Let $\mathcal{ M = M}^{(\Phi , \Psi)}$ and $\mathcal{O = O}^{(\Phi , \Psi)}$. Let $O \in \BL(\Hil_1,\Hil_2)$, then

$$\left( \mathcal{M} \left( O \right) c\right)_j = \sum \limits_k \left( \mathcal{M} \left( O \right) \right)_{j,k} c_k =  \sum \limits_k \left<O \psi_k, \phi_j \right> c_k = $$
\begin{equation}\label{sec:basicproof1} =  \left< \sum \limits_k  c_k O \psi_k, \phi_j \right> = \left< O \sum \limits_k c_k \psi_k, \phi_j \right> = \left< O D_{\Psi} c, \phi_j \right>
\end{equation}
$$ \Longrightarrow \norm{2}{\mathcal M c}^2 = \sum \limits_j \left| \left< O D_{\Psi} c, \phi_j \right> \right|^2 \le  B' \cdot \norm{\Hil}{ O D_{\Psi} c}^2 \le B' \cdot \norm{Op}{O}^2 B \norm{2}{c}^2 $$
Equation \ref{sec:basicproof1} also shows us, that as operator we have 
$$ {\mathcal M}^{(\Phi , \Psi)} \left( O \right) = C_{\Phi} \circ O \circ D_{\Psi} .$$

On the other hand let $M$ be an infinite matrix, then
$$ \mathcal{O} \left( M \right) = D_{\Phi} \circ \mathcal M \circ C_{\Psi}$$ 
$$ \Longrightarrow \norm{\Hil_1 \rightarrow \Hil_2}{ \mathcal{O} \left( M \right)} \le \norm{l^2 \rightarrow \Hil_2}{D_{\Phi}} \cdot \norm{l^2 \rightarrow l^2}{M} \cdot \norm{\Hil_1 \rightarrow l^2}{C_{\Psi}} \le $$
$$ \le \sqrt{B'} \cdot \norm{l^2 \rightarrow l^2}{M} \sqrt{B}$$
\end{proof}

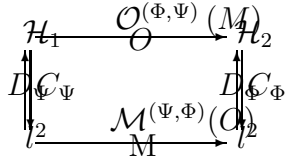
\begin{figure}[ht]
\center
	\begin{picture}(100,60)
	\put(10,50){$\Hil_1$}
	\put(90,50){$\Hil_2$}
	\put(15,52){\vector(4,0){72}}
	\put(45,55){$\mathcal{O}^{(\Phi , \Psi)} \left( M \right)$} 
	\put(50,47){$O$} 
	\put(13,47){\vector(0,-4){30}}
	\put(11,17){\vector(0,4){30}}
	\put(93,47){\vector(0,-4){30}}
	\put(91,17){\vector(0,4){30}}
	\put(11,10){$l^2$}
	\put(4,30){$D_{\Psi}$}
	\put(15,30){$C_{\Psi}$}
	\put(84,30){$D_{\Phi}$}
	\put(95,30){$C_{\Phi}$}
	\put(91,10){$l^2$}
	\put(15,12){\vector(4,0){72}}
	\put(50,7){M} 
	\put(43,17){$\mathcal M^{(\Psi,\Phi)}(O)$}
	\end{picture}
\caption{The operator induced by a matrix $M$ and the matrix induced by an operator $O$} \label{fig:matop1}
\end{figure}

\begin{definition} For an operator $O$ and a matrix $M$ as in Theorem \ref{sec:matbyfram1}, we call ${\mathcal M}^{(\Psi,\Phi)}(O)$ the \em matrix induced by the operator $O$ \em with respect to the Bessel sequences $\Psi = (\psi_k)$ and $\Phi = (\phi_k)$ 
and ${\mathcal O}^{(\Psi,\Phi)}(M)$ the \em operator induced by the matrix $M$ \em with respect to the Bessel sequences $\Psi$ and $\Phi$. 
(See Figure \ref{fig:matop1}.)
\end{definition} 

If we do not want to stress the dependency on the frames and there is no change of confusion, the notation ${\mathcal M} (O)$ and ${\mathcal O} (M)$ will be used.

For frames we can prove more properties: 
\begin{proposition} \label{sec:propmatropfram1} Let $\Psi = (\psi_k)$ be a frame in $\Hil_1$ with bounds $A,B$, $\Phi = (\phi_k)$ in $\Hil_2$ with $A',B'$. Then
\begin{enumerate}
\item $ \left( {\mathcal O^{(\Phi , \Psi)} \circ M^{(\tilde{\Phi}, \tilde{\Psi})}}\right)  = Id = \left( {\mathcal O^{(\tilde{\Phi}, \tilde{\Psi})} \circ M^{(\Phi , \Psi)}}\right) $. \\
And therefore for all $O \in \BL(\Hil_1,\Hil_2)$:
$$ O = \sum \limits_{k,j} \left<O \tilde{\psi}_j, \tilde{\phi}_k \right>  \phi_k \otimes_i \overline{\psi}_j $$
\item  $\mathcal{M}^{(\Phi , \Psi)}$ is injective and $\mathcal{O}^{(\Phi , \Psi)}$ is surjective.
\item Let $\Hil_1 = \Hil_2$, then $\mathcal{O}^{(\Psi, \tilde{\Psi})} (Id_{l^2}) = Id_{\Hil_1}$ 
\item Let $\Xi = (\xi_k)$ be any frame in $\Hil_3$, and $O : \Hil_3 \rightarrow \Hil_2$ and $P: \Hil_1 \rightarrow \Hil_3$. Then
$$ \mathcal{M}^{(\Phi, \Psi)}\left( O \circ P \right) = \left( \mathcal{M}^{(\Phi, \Xi)}\left( O \right) \cdot \mathcal{M}^{(\tilde{\Xi}, \Psi)} \left( P \right) \right) $$
\end{enumerate}

\end{proposition}
\begin{proof} 
1.) For $f \in \Hil_1$ we get  $ \left( \mathcal O^{(\Phi , \Psi)} \circ M^{(\tilde{\Phi}, \tilde{\Psi})} \right) \left(O\right) (f) =  \sum \limits_k  \left( \sum \limits_j \left<O \tilde{\psi}_j, \tilde{\phi}_k \right> \left<f, \psi_j\right> \right) \phi_k  = $
$$ = \sum \limits_k  \left( \left< \sum \limits_j \left<f, \psi_j\right> O \tilde{\psi}_j, \tilde{\phi}_k \right>  \right) \phi_k = 
\sum \limits_k  \left< O f , \tilde{\phi}_k \right> \phi_k = Of .$$

For the other equality the roles of the frame and the dual just have to be switched.

2.) From $\mathcal{O M} = Id$ we know that $\mathcal{M}$ is injective and $\mathcal{O}$ is surjective.

3.) $\mathcal{O} (Id) f = \sum \limits_k  \left( \sum \limits_j \delta_{k,j} \left<f, \tilde{\psi}_j\right> \right) \psi_k  = \sum \limits_k  \left<f, \tilde{\psi}_k \right> \psi_k = f $.

4.)  $  \mathcal{M}^{(\Phi, \Psi)}\left( O \circ P \right)_{p,q} = \left< O \circ P \psi_q, \phi_p \right> =  \left< P \psi_q, O^* \tilde{\phi}_p \right> $. 

On the other hand
$$ \left( \mathcal{M}^{(\Phi, \Xi)}\left( O \right) \cdot \mathcal{M}^{(\tilde{\Xi}, \Psi)} \left( P \right) \right)_{p,q} = \sum \limits_k \mathcal{M}^{(\Phi, \Xi)} \left( O \right)_{p,k} \cdot \mathcal{M}^{(\tilde{\Xi}, \Psi)}  \left( P \right)_{k,q} = $$
$$ = \sum \limits_k \left< O h_k, \phi_p\right> \left< P \psi_q, \tilde{\xi}_k \right> = \sum \limits_k \left< h_k, O^* \tilde{\phi}_p\right> \left< P \psi_q, \tilde{\xi}_k \right> = $$
$$ = \left<  \sum \limits_k \left< P \psi_q, \tilde{\xi}_k \right> \xi_k, O^* \tilde{\phi}_p \right> = \left< P \psi_p , O^* \tilde{\phi}_p \right> .$$
\end{proof}

As a direct consequence we get the following corollary:
\begin{corollary} For the frame $\Phi = (\phi_k)$ the function $\mathcal M^{(\Phi,\tilde{\Phi})}$ is a Banach-algebra monomorphism between the algebra of bounded operators $\left( \BL(\Hil_1,\Hil_1), \circ\right)$ and the infinite matrices of $\left( \BL(l^2,l^2), \cdot \right)$.  
\end{corollary}

The other function $\mathcal{O}$ is in general not so ``well-behaved''. It is, if the dual frames are biorthogonal. In this case these functions are isomorphisms, refer to Section \ref{sec:matrepriesz0}. 

\begin{lemma} \label{sec:matropsurinj1} Let $O : \Hil_1 \rightarrow \Hil_2$ be a linear and bounded operator, let $\Psi = (\psi_k)$ and $\Phi = (\phi_k)$ be frames in $\Hil_1$ resp. $\Hil_2$. Then $\mathcal M^{(\Phi , \tilde{\Psi})}(O)$ maps $ran\left( C_{\Psi}\right)$ into $ran \left( C_{\Phi} \right)$ with
$$ \left( \left< f , \psi_k \right> \right)_k \mapsto \left( \left< O f , \phi_k \right> \right)_k .$$
If $O$ is surjective, then $\mathcal M^{(\Phi , \tilde{\Psi})}  (O)$ maps $ran\left( C_{\Psi}\right)$ onto $ran \left( C_{\Phi} \right)$.
If $O$ is injective, $\mathcal M^{(\Phi , \tilde{\Psi})}  (O)$ is also injective.
\end{lemma} 
\begin{proof}
Let $c \in ran(C_{\Psi})$, then there exists $f \in \Hil_1$ such that $c_k = \left< f, \psi_k \right>$. 
$$\left( \mathcal{M}^{(\Phi, \tilde{\Psi})} (O) (c) \right)_i = \sum \limits_k \left< O \tilde{\psi}_k, \phi_i \right>\left< f, \psi_k \right> 
= \left< \sum \limits_k  \left< f, \tilde{\psi}_k \right>  O \psi_k, \phi_i \right> = \left< O f , \phi_i \right> $$
So $ \left( \left< f, \psi_k \right> \right)_{k} \mapsto \left( \left< O f , \phi_k \right> \right)_k$.

If $O$ is surjective, then for every $f$ there exists a $g$ such that $Og = f$, and therefore $ \left< g, \psi_k \right> \mapsto \left< f , \phi_k \right>$.

If $O$ is injective, then let's suppose that $ \left< O f , \phi_k \right> =  \left<O g , \phi_k \right>$. Because $(\phi_k)$ is a frame $ \Longrightarrow O f =  O g$ $\Longrightarrow f =  g$ $\Longrightarrow \left< f, \psi_k \right> = \left< g, \psi_k \right>$.
\end{proof}

Particularly for $O = Id$ the Gram matrix $G_{\Phi,\tilde{\Psi}} = \left( \left< \tilde{\psi}_k , \phi_i \right> \right)_{k,i}$ maps $ran\left( C_{\Psi}\right)$ bijectively on $ran\left( C_{\Phi}\right)$. So we get a way to a way to ``switch'' between frames by mapping from one analysis range into the other \cite{xxlphd1}.

Let us give some examples:
\Bsp{sec:bspopfram1}  Let $\Psi = (\psi_k)$ and $\Phi = (\phi_k)$ be frames in $\Hil$ 
and $\Delta = (\delta_k)$ the canonical basis of $l^2$. Then
\begin{enumerate}
\item $S_{\Psi}: \Hil \rightarrow \Hil$ $\Longrightarrow$ ${\mathcal M^{(\Psi,\tilde{\Psi})} (S_{\Psi})} = G_{\Psi}$.
\item $S_{\Psi}^{-1} : \Hil \rightarrow \Hil$ $\Longrightarrow$ ${\mathcal M^{(\Psi,\tilde{\Psi})} (S_{\Psi}^{-1})} = G_{\tilde{\Psi}}$.
\item $C_{\Phi}: \Hil \rightarrow l^2$ $\Longrightarrow$
$$\mathcal M^{(\Delta, \Psi)} (C_{\Phi})_{k,j} = \left< C_{\Phi} \psi_j, \delta_k \right> = \sum \limits_l \left<\psi_j , \phi_l \right> \left< \delta_{k}, \delta_l\right> = \left<\psi_j , \phi_k \right> = \left( G_{\Phi, \Psi}\right)_{k,j}.$$
\item $Id : \Hil \rightarrow \Hil$ $\Longrightarrow$ $\mathcal M^{(\Phi,\Psi)} (Id) = G_{\Phi, \Psi}$.
\item $Id : l^2 \rightarrow l^2$ $\Longrightarrow$ $\mathcal O^{(\Phi,\Psi)} (Id) = D_{\Phi} \circ C_{\tilde{\Psi}} = S_{\Psi, \tilde{\Phi}}$.
\end{enumerate}

\subsubsection{Motivation: Solving Operator Equalities} \label{sec:soopeq0}

Given an operator equality 
\begin{equation} \label{sec:solvopeq1}
~ O \cdot f = g ~
\end{equation}
it is natural to discretize it to find a solution. Let $\Phi = (\phi_k)$ be a frame. Let us suppose that for a given $g$ with coefficients $d = (d_k) = ( \left< g , \phi_k \right> )$ and a matrix representation $M$ of $O$ there is an algorithm 
to find the least square solution of 
\begin{equation} \label{sec:solvopmateq1}
~ M \cdot c = d ~
\end{equation}
for example using the pseudoinverse \cite{ole1}. 
Still, if using frames, 
we can not expect to find a true solution for Eq. \ref{sec:solvopeq1} just by applying $D_{\tilde \Phi}$ on $c$ as in general $c$ is not in $ran (C_{\Phi})$ even if $d$ is. But rephrasing Eq. \ref{sec:solvopeq1} we see the following: 
$$ O f = g \Longleftrightarrow \sum \limits_k \left< f, \phi_k\right> O \tilde \phi_k = g  \Longleftrightarrow \sum \limits_k \left< f, \phi_k\right> \left< O \tilde \phi_k , \phi_k \right> = \left< g , \phi_k \right>  $$
$$ \Longleftrightarrow {\mathcal M}^{(\Phi , \tilde{\Phi})} \left( O \right) \cdot C_\Phi f = C_\Phi g.$$
It can be easily seen that this is equivalent to projecting $c$ on $ran(C)$, solving $M C_{\Phi} D_{\tilde{\Phi}} c =d$, which is a common idea found in many algorithms, for example for a recent one see \cite{teschk1}.
\\

This gives us an algorithm for finding an approximative solution to the inverse operator problem $O f = g$. 
\begin{enumerate}
\item Set $M = {\mathcal M}^{(\Phi , \tilde{\Phi})} \left( O \right)$. 
\item Find a good finite dimensional approximation $M_N$ of $M$ by using the finite section method \cite{gohberg1} and 
\item then apply an algorithm like e.g. the QR factorization \cite{trebau1} to find a solution for Eq. \ref{sec:solvopeq1}. 
\item 
and synthezise with the dual frame $\tilde{\Phi}$.
\end{enumerate}

\Rem It has been shown in \cite{olestroh1}, that the finite section is very useful in the case of frame theory. It would be very interesting to investigate the idea presented above further in this context.

\subsection{Matrix representation using Riesz Bases} \label{sec:matrepriesz0}
The coefficients using a Riesz basis are unique, so Theorem \ref{sec:matbyfram1} can be extended to:
\begin{theorem} \label{sec:rieszmatrixrep1} 
 Let $\Phi = (\phi_k)$ be a Riesz basis for $\Hil_1$, $\Psi = (\psi_k)$ one  for $\Hil_2$. The functions $\mathcal{M}^{(\Phi,\Psi)}$ and  $\mathcal{O}^{(\tilde{\Phi},\tilde{\Psi})}$ between 
$\BL(\Hil_1,\Hil_2)$ and the infinite matrices in $\BL(l^2,l^2)$
 are bijective. $\mathcal{M}^{(\Phi,\Psi)}$ and $\mathcal{O}^{(\tilde{\Phi},\tilde{\Psi})}$ are inverse to each other. For $\Hil_1 = \Hil_2$ 
 the identity is mapped on the identity by $\mathcal{M}^{(\Phi,\Psi)}$ and $\mathcal{O}^{(\tilde{\Phi},\tilde{\Psi})}$. If furthermore $\Psi = \Phi$ then 
$\mathcal{M}^{(\Phi,\tilde{\Phi})}$ and $\mathcal{O}^{(\Phi,\tilde{\Phi})}$ are Banach algebra isomorphisms, respecting the identities $id_{l^2}$ and $id_\Hil$.
\end{theorem}
\begin{proof} We know from Proposition \ref{sec:propmatropfram1} that ${\mathcal O \circ \mathcal M } = Id$. Let's look at
$$ {\left( \left( {\mathcal M \circ \mathcal O} \right) \left(M \right)\right)}_{p,q} =  {\mathcal M}\left( \sum \limits_k \sum \limits_j M_{k,j} \left< \cdot , \psi_j\right> \phi_k \right)_{p,q} = $$
$$ = \left< \sum \limits_k \sum \limits_j M_{k,j} \left< \tilde{\psi}_q, \psi_j\right> \phi_k , \tilde{\phi}_p \right>  =  \sum \limits_k \sum \limits_j  M_{k,j}   \underbrace{\left< \tilde{\psi}_q, \psi_j\right>}_{\delta_{k,p}} \underbrace{\left< \phi_k , \tilde{\phi}_p \right>}_{\delta_{k,p}} = M_{p,q} $$
So these functions are inverse to each other and therefore bijective.

$$ {\mathcal M}(Id_{\Hil \rightarrow \Hil})_{p,q} =  \left< \psi_q, \tilde{\psi}_p \right> = \delta_{q,p}  = \left( Id_{l^2 \rightarrow l^2} \right)_{p,q} $$

We know that $\mathcal{M}^{(\Phi,\tilde{\Phi})}$ is a Banach algebra homomorphism and so is its inverse.
\end{proof}

\subsection{Matrix Representation of $\HS$ Operators} \label{sec:descrhsopfram0} 

We now have the adequate tools to state that $\HS$ operators correspond exactly to the Frobenius matrices, as expected. 
\begin{definition}\label{sec:frobnorm1}
Let $A$ be an $m$ by $n$ matrix, then
$$ \left\| A \right\|_{fro} = \sqrt{\sum \limits_{i=0}^{n-1} \sum \limits_{j=0}^{m-1} \left| a_{i,j} \right|^2} $$
is the \em Frobenius 
norm\em . 
Let us denote the set of all matrices with finite Frobenius norm by $l^{(2,2)}$, the set of \em Frobenius matrices\em .
\end{definition}

\begin{proposition} \label{sec:frobmatrHSop1} Let $\Psi = (\psi_k)$ be a Bessel sequence in $\Hil_1$ with bound $B$, $\Phi = (\phi_k)$ in $\Hil_2$ with $B'$. Let $M$ be a matrix in $l^{(2,2)}$. Then $\mathcal{O}^{(\Phi,\Psi)}(M) \in \HS(\Hil_1, \Hil_2)$, the Hilbert Schmidt class of operators from $\Hil_1$ to $\Hil_2$, with $\norm{\HS}{ \mathcal{O}(M) } \le \sqrt{B B'} \norm{fro}{M}$.

Let $O \in \HS$, then $\mathcal{M}^{(\Phi,\Psi)}(O) \in l^{(2,2)}$ with $\norm{fro}{\mathcal{M}(O)} \le \sqrt{B B'} \norm{\HS}{O}$.
\end{proposition}
\begin{proof} 1.) Naturally the matrices in $l^{(2,2)}$ correspond to Hilbert-Schmidt operators on $l^2$ as
$ \left\| M \right\|_{HS(l^2)}^2 = \sum \limits_i \left\| M \delta_i \right\|_{\Hil_1}^2 = \sum \limits_i \sum \limits_p \left| M_{p,i} \right|^2 =  \left\| M \right\|_{fro}$.
As the Hilbert-Schmidt class of operators is an ideal, we know that
$$ \norm{\HS(\Hil_1,\Hil_2)}{{\mathcal O}(M)} = \norm{\HS(\Hil_1,\Hil_2)}{D_{\Phi} \circ M \circ C_{\Psi}} \le $$
$$ \le \norm{Op(l^2,\Hil_2)}{D_{\Phi}} \cdot \norm{\HS(l^2)}{ M } \cdot \norm{Op(\Hil_1,l^2)}{C_{\Psi}} = \sqrt{B B'} \norm{fro}{M} $$

2.) Let $(e_p)$ be any ONB of $\Hil_2$.
$$ \norm{fro}{\mathcal{M}(O)}^2 = \sum \limits_l \sum \limits_k \left| \left< O \psi_l, \phi_k \right> \right|^2 \le \sum \limits_l B' \cdot \norm{\Hil_1}{O \psi_l}^2 = \sum \limits_l B' \sum \limits_p \left| \left< O \psi_l, e_p\right> \right|^2 = $$
$$ = \sum \limits_p B' \sum \limits_l \left| \left< \psi_l, O^* e_p\right> \right|^2 \le \sum \limits_p B' B \norm{\Hil_2}{O^* e_p}^2 \le B' B \norm{\HS}{O^*}^2 =  B B' \norm{\HS}{O}^2 $$
\end{proof}

\subsubsection{Matrices and the Kernel Theorems}

For $L^2(\Rd)$ the $\HS$ operators are exactly those integral operators with kernels in $\LtRtd$ \cite{feikoz1,schatt1}. This means that there exists a $\kappa_O \in L^2(\RR^{2d})$ such an operator can be described as
$$ \left( O f \right) (x) = \int \kappa_O(x,y) f(y) dy$$
Or in weak formulation
\begin{equation} \label{sec:weakkernel1}
\left< O f , g \right> = \int \int \kappa_O(x,y) f(y) \overline{g}(x) dy dx = \left< \kappa_O , f \otimes_o \overline{g} \right>.
\end{equation}
From \ref{sec:matbyfram1} we know that 
$$ O = \sum \limits_{j,k} \left<O \tilde{\psi}_j, \tilde{\phi}_k \right>  \phi_k \otimes_i \overline{\psi}_j $$
and so
\begin{corollary} Let $O \in \HS\left( L^2 \left( \Rd \right) \right)$. Let $\Psi = (\psi_j)$ and $\Phi = (\phi_k)$ be frames in $\LtRd$. Then the kernel of $O$ is given as:
$$ \kappa_O = \sum \limits_{j,k} \left<O \tilde{\psi}_j, \tilde{\phi}_k \right> \cdot \phi_k \otimes_o \overline{\psi}_j = \sum \limits_{j,k} {\mathcal M}^{(\tilde{\Psi},\tilde{\Phi})}(O)_{k,j} \cdot \phi_k \otimes_o \overline{\psi}_j $$
\end{corollary}
\begin{proof}
$$ \kappa(O) = \kappa \left( \sum \limits_{j,k} \left<O \tilde{\psi}_j, \tilde{\phi}_k \right>  \phi_k \otimes_i \overline{\psi}_j  \right) = $$
$$ = \sum \limits_{j,k} \left<O \tilde{\psi}_j, \tilde{\phi}_k \right> \kappa \left(  \phi_k \otimes_i \overline{\psi}_j  \right) =  \sum \limits_{j,k} \left<O \tilde{\psi}_j, \tilde{\phi}_k \right> \phi_k \otimes_o \overline{\psi}_j $$
as  $\left( f \otimes_i g \right) (h) = \left< h , g \right> f = \int h(x) \overline{g(x)} dx \dot f(y) $ and so $\kappa (f \otimes_i g) = f \otimes_o g$. 

\end{proof}

\section{Perspectives} \label{sec:concpersp0}

Apart from $\HS$ operators on $L^2(\Rd)$, there is a large variety of function spaces, where operators are exactly integral operator using the weak formulation in Eq. (\ref{sec:weakkernel1}), for example for operators on the Schwartz space $O: \mathcal{S \rightarrow S'}$ , on  modulation spaces $O: M_v^1(\Rd) \rightarrow M_{1/v}^\infty$ \cite{Groech1},or on Feichtinger's algebra $O: S_0' \rightarrow S_0$ and $O: S_0 \rightarrow S_0'$ \cite{feikoz1}. 
In order to derive similar results for the case of Banach spaces of functions or distributions, Section \ref{sec:descropfram0} could be generalized to these function spaces.

Using the connection between operators and frames, we can ask, which operators are induced by diagonal matrices. Let $m$ be a sequence and $\mathfrak{diag}(m)$ the matrix that has this sequence as diagonal. Then define
$$ \MM_{m,\Phi,\Psi} := \mathcal{O}^{(\Phi , \Psi)} \left( \mathfrak{diag}(m) \right) = \sum \limits_k m_k \cdot  \phi_k \otimes \psi_k$$
This means we have arrived quite naturally at the definition of frame multiplier as introduced in \cite{xxlmult1}. This connection should be investigated and exploited further.

Especially in areas where orthonormal bases are not very useful and frames are already used heavily for analysis-synthesis systems
, the idea of matrix representation with frames should be used in applications. This concept should be applied using the finite section method \cite{gohberg1} as mentioned in Section \ref{sec:soopeq0}. Furthermore, to conserve the semantic connection between model and discretization even more, it might be interesting to use frames in a projection method like the finite section method for approximation of infinite matrices by finite ones. Instead of using the ONB for the projection use the frames, which are already connected to the problem. 
In this sense the ideas in \cite{olestroh1} could be extended to the notion of self-localized frames as presented in \cite{forngroech1}.

\vspace{13pt}
\section{ACKNOWLEDGEMENT}
\vspace{13pt}
The author would like to thank Hans G. Feichtinger, Bruno Torr\'esani and Jean-Pierre Antoine for
many helpful comments and suggestions as well as Wolfgang Kreuzer for proofreading. He would like to thank the hospitality of the LATP, CMI, Marseille, France and FYMA, UCL, Louvain-la-Neuve, Belgium, where part of this work was prepared, supported by the HASSIP-network. 
This work was partly supported by the European Union's Human Potential Program, under contract HPRN-CT-2002-00285 (HASSIP).

\small

\end{document}